**An explicit mapping between the frequency domain and the time domain representations of nonlinear systems**

Marissa Condon[1,*] and Rossen Ivanov[2,3]

**[1]School of Electronic Engineering, Dublin City University, Dublin 9, IRELAND**

**[2]School of Mathematics, Trinity College, Dublin 2, IRELAND**

**[3]Institute for Nuclear Research and Nuclear Energy, 72 Tzarigradsko chaussee, 1784 Sofia, BULGARIA**

[*]Correspondence to: Marissa Condon, School of Electronic Engineering, Dublin City University, Dublin 9, Ireland; Phone: +353 1 700 5405, Fax: +353 1 700 5508, E-mail: marissa.condon@dcu.ie

e-submission number: **NLM_rossen.ivanov.AT.dcu.ie_20050119/1**

ABSTRACT

Explicit expressions are presented that describe the input-output behaviour of a nonlinear system in both the frequency and the time domain. The expressions are based on a set of coefficients that do not depend on the input to the system and are universal for a given system. The anharmonic oscillator is chosen as an example and is discussed for different choices of its physical parameters. It is shown that the typical approach for the determination of the Volterra Series representation is not valid for the important case when the nonlinear system exhibits oscillatory behaviour and the input has a pole at the origin (in the frequency domain), e.g. the unit-step function. For this case, resonant effects arise and the analysis requires additional care.

# 1. Introduction

The mathematical analysis of linear systems in both the frequency and the time domain is a well-developed area but nevertheless is still under active investigation. Various relationships between the frequency and time domain descriptions are constantly developed and reformulated for applications in different situations e.g.[1-5] to mention only a few.

The situation is much more complex for nonlinear systems. Most input-output representations for nonlinear systems are based on the Volterra series expansion [6]. The transformation between the frequency and time domain is not, in general, straightforward [7]. The Volterra kernels arising from multidimensional Fourier transforms are often related to the so-called generalized frequency response functions and desribing functions [8,9]. In [10-12], an analytical relationship between non-linear integro-differential equations and the generalized frequency response functions is derived and further analysed.

In this work an alternative description is suggested that is based on the Fourier coefficients (henceforth termed 'scattering' coefficients) arising from the expansion of the Volterra kernels. The 'scattering' coefficients are universal for each nonlinear system and provide a unique description of the input-output behaviour of the system in both the frequency and time domains. Also, if the nonlinear equation governing the system is explicitly given, these coefficients can easily be computed via a simple numerical algorithm.

In Section 2 an explicit mapping between the frequency domain and the time domain is presented. The computation of the 'scattering' coefficients for the mapping from the Volterra kernels associated with a given differential equation is given in Section 3.

The example of the anharmonic oscillator is discussed in Section 4. A numerical algorithm for the computation of the 'scattering' coefficients directly from the differential equation governing the system is outlined in Section 5. Section 6 addresses the inverse problem related to the system identification: how to compute the 'scattering' coefficients from input-output measurements when the differential equation describing the (in general-nonlinear) system is unknown. This is of course a challenging question that requires additional investigation and may necessitate making certain assumptions regarding the nonlinear system in question. The inverse problem is addressed in this work in relation to linear systems and is presented as a starting point for future work.

## 2. The mapping between the frequency domain and the time domain for nonlinear systems.

Consider a nonlinear system in the frequency domain with an input and output $X(\omega)$ and $Y(\omega)$ correspondingly. Let the input-output mapping be given by the Volterra Series Expansion:

$$Y(\omega) = H_1(\omega)X(\omega) + \frac{1}{2\pi}\int H_2(\omega_1,\omega_2)X(\omega_1)X(\omega_2)\delta(\omega_1+\omega_2-\omega)d\omega_1 d\omega_2 + \ldots$$
$$+ \frac{1}{(2\pi)^{n-1}}\int H_n(\omega_1,\ldots,\omega_n)X(\omega_1)\ldots X(\omega_n)\delta(\omega_1+\ldots+\omega_n-\omega)d\omega_1\ldots d\omega_n + \ldots \quad (1)$$

where $H_n(\omega_1,\ldots,\omega_n)$ is a kernel of degree $n$ and a completely symmetric function of its arguments. All integrations are from $-\infty$ to $\infty$.

Suppose that $H_n(\omega_1,\ldots,\omega_n)$ are nonzero when their arguments take values in the interval $[-\omega_M,\omega_M]$ where $\omega_M$ is assumed to be large, but finite. Then $H_n(\omega_1,\ldots,\omega_n)$ may be expanded in a Fourier series as follows:

$$H_n(\omega_1,\ldots,\omega_n) = \sum_{k_1,\ldots,k_n} a^{[n]}_{k_1,\ldots,k_n} e^{-jk_1\tilde{\omega}_1-\ldots-jk_n\tilde{\omega}_n} \qquad (2)$$

where all summation indices run from zero to infinity, $a^{[n]}_{k_1,\ldots,k_n}$ is real and completely symmetric with respect to its lower indices and $\tilde{\omega} = \pi\omega/\omega_M$. The representation in (2) is valid only within the interval $[-\omega_M, \omega_M]$ since as already stated it is assumed that $H_n(\omega_1,\ldots,\omega_n) = 0$ otherwise.

Suppose that an (arbitrary) input $x(t)$ is applied at $t = 0$ (i.e. input signal $x(t)\theta(t)$ where $\theta(t)$ is the unit step-function). In what follows, for simplicity, the Fourier transform of $x(t)\theta(t)$ will be denoted by $X(\omega)$. Then the output $y(t)$ is given by the following:

**Proposition 1.**

$$y(t) = \sum_{n=1}^{\infty} \sum_{k_1,\ldots,k_n} a^{[n]}_{k_1,\ldots,k_n} \prod_{i=1}^{n} x(t-\tilde{k}_i)\theta(t-\tilde{k}_i) \qquad (3)$$

where $\tilde{k} = \dfrac{\pi k}{\omega_M}$.

*Proof*: Taking (1) and (2) into account and using the following representation of the delta-function:

$$\delta(\omega) = \frac{1}{2\pi}\int_{-\infty}^{\infty} e^{j\omega v} dv \qquad (4)$$

it follows that

$$y(t) = \frac{1}{2\pi}\int_{-\infty}^{\infty} e^{j\omega t} Y(\omega) d\omega$$

$$= \frac{1}{2\pi}\sum_{k=0}^{\infty} a^{[1]}_k \int_{-\infty}^{\infty} e^{j\omega(t-\tilde{k})} X(\omega)d\omega +$$

$$\sum_{n=2}^{\infty}\frac{1}{(2\pi)^n}\sum_{k_1,\ldots,k_n} a^{[n]}_{k_1,\ldots,k_n}\int_{-\infty}^{\infty}\ldots\int_{-\infty}^{\infty} e^{j\omega t}\left(\frac{1}{2\pi}\int_{-\infty}^{\infty} e^{j(\omega_1+\ldots+\omega_n-\omega)v}dv\right)\prod_{i=1}^{n}\left(e^{-j\omega_i\tilde{k}_i} X(\omega_i)d\omega_i\right)d\omega$$

$$= \sum_{k=0}^{\infty} a_k^{[1]} x(t-\tilde{k})\theta(t-\tilde{k}) +$$

$$\sum_{n=2}^{\infty} \frac{1}{(2\pi)^n} \sum_{k_1,\ldots,k_n} a_{k_1,\ldots,k_n}^{[n]} \int_{-\infty}^{\infty} \cdots \int_{-\infty}^{\infty} \left( \frac{1}{2\pi} \int_{-\infty}^{\infty} e^{j\omega(t-v)} d\omega \right) \prod_{i=1}^{n} \left( e^{j\omega_i(v-\tilde{k}_i)} X(\omega_i) d\omega_i \right) dv$$

$$= \sum_{k=0}^{\infty} a_k^{[1]} x(t-\tilde{k})\theta(t-\tilde{k}) +$$

$$\sum_{n=2}^{\infty} \sum_{k_1,\ldots,k_n} a_{k_1,\ldots,k_n}^{[n]} \int_{-\infty}^{\infty} \delta(t-v) \prod_{i=1}^{n} \left( \frac{1}{2\pi} \int_{-\infty}^{\infty} e^{j\omega_i(v-\tilde{k}_i)} X(\omega_i) d\omega_i \right) dv$$

$$= \sum_{n=1}^{\infty} \sum_{k_1,\ldots,k_n} a_{k_1,\ldots,k_n}^{[n]} \prod_{i=1}^{n} x(t-\tilde{k}_i)\theta(t-\tilde{k}_i).$$

**Corollary 1**. The representation in (2) represents a causal input-output mapping due to (3).

**Definition 1**. The set of coefficients $a_{k_1,\ldots,k_n}^{[n]}$ are called 'scattering' coefficients.

**Corollary 2.** The 'scattering' coefficients uniquely determine the input-output mapping of the nonlinear system. They provide an explicit parametrization for this mapping in both the frequency (1)-(2) and the time domain (3).

**Corollary 3**. For linear systems $H_1(\omega)$ is simply the transfer function ($a_k \equiv a_k^{[1]}$):

$$H_1(\omega) = \sum_{k=0}^{\infty} a_k e^{-jk\tilde{\omega}} \tag{5}$$

for $\omega \in [-\omega_M, \omega_M]$ and zero everywhere else.

$$\operatorname{Re} H_1(\omega) = \sum_{k=0}^{\infty} a_k \cos k\tilde{\omega} \tag{6}$$

$$\operatorname{Im} H_1(\omega) = -\sum_{k=0}^{\infty} a_k \sin k\tilde{\omega} \tag{7}$$

The expression in (6) describes an even function and the expression in (7) describes an odd function. Both are defined for $\omega \in [-\omega_M, \omega_M]$, (i.e. $\tilde{\omega} \in [-\pi, \pi]$) and zero everywhere else. Due to causality, the expressions in (6) and (7) are linked via the Kramers-Kronig relations (Hilbert transform) [4, 5, 13, 14].

Therefore, once the 'scattering' coefficients are obtained (e.g. from measurements in the frequency domain (1)), then the response for an arbitrary input may easily be determined from (3). The meaning of the Fourier coefficients $a^{[n]}_{k_1,\ldots,k_n}$ in the time-domain is clearly evident from (3) - they give the contribution of the input at a given time, preceding the output, to the output signal.

### 3. Computation of the 'scattering' coefficients - direct problem

Suppose that the nonlinear differential equation that describes the system is explicitly given. Then, then the computation of the 'scattering' coefficients may be implemented as follows.

First, the kernels $H_n(\omega_1,\ldots,\omega_n)$ are computed following a procedure such as that outlined in [10]. Then, once the kernels are determined, the 'scattering' coefficients may be evaluated as follows: Consider, for example, the degree one kernel given in (5). One can write:

$$\int_{-\infty}^{\infty} e^{jl\tilde{\omega}} H_1(\omega) d\omega = \int_{-\omega_M}^{\omega_M} e^{jl\tilde{\omega}} H_1(\omega) d\omega =$$

$$= 2\omega_M a_l + \sum_{k \neq l} a_k \int_{-\omega_M}^{\omega_M} e^{j(\tilde{l}-\tilde{k})\omega} d\omega$$

$$= 2\omega_M a_l + 2\sum_{k \neq l} a_k \frac{\sin((l-k)\pi)}{\tilde{l}-\tilde{k}}$$

$$= 2\omega_M a_l$$

and therefore:

$$a^{[1]}_l \equiv a_l = \frac{1}{2\omega_M} \int_{-\infty}^{\infty} e^{jl\tilde{\omega}} H_1(\omega) d\omega \qquad (8)$$

Similarly, the above derivation can easily be generalised to enable the determination of the higher order scattering coefficients:

$$a^{[n]}_{k_1\ldots k_n} = \frac{1}{(2\omega_M)^n} \int_{-\infty}^{\infty}\ldots\int_{-\infty}^{\infty} e^{jk_1\tilde{\omega}_1+\ldots+jk_n\tilde{\omega}_n} H_n(\omega_1,\ldots,\omega_n) d\omega_1\ldots d\omega_n \quad (9)$$

The expressions in (8) and (9) are convenient for computation with the Residue Theorem. Note that in reality the kernels are not identically zero outside a finite interval. However, it is always the case that it is possible to identify a finite interval, i.e. a bound $\omega_M$ ('the highest frequency in the system'), such that the contribution of the kernels $H_n(\omega_1,\ldots,\omega_n)$ for frequencies greater than $\omega_M$ in (9) is negligible. Obviously, with increasing $\omega_M$, the precision in the computation of the 'scattering' coefficients is increased. Note that in order to cover a fixed time interval $[0,T_M]$ in the time domain (3), a finite number of 'scattering' coefficients is required, namely $a^{[n]}_{k_1,\ldots,k_n}$ with $k_i \leq N$, where $\tilde{N} = T_M$, i.e.

$$N = T_M \omega_M / \pi \quad (10)$$

Therefore, the number of coefficients $a^{[n]}_{k_1,\ldots,k_n}$ necessary for a simulation within a fixed and finite time interval $T_M$, depends on $\omega_M$. Increasing $\omega_M$ increases the precision of the simulation, but more coefficients are required according to (10).

### 4. Examples

As an example, consider the very well known anharmonic oscillator. The nonlinear equation governing the behaviour of anharmonic oscillator with damping is:

$$\ddot{y} + b\dot{y} + \omega_0^2 y + \varepsilon y^2 = x(t)\theta(t) \quad (11)$$

where $b > 0$ is the damping constant, $\varepsilon$ is the anharmonicity coefficient that is supposed to be small enough such that the perturbative expansion in (1) is valid; $\omega_0$ is the frequency of the oscillator when $b = \varepsilon = 0$.

The following cases are possible: an underdamped oscillator ($\omega_0 > b/2$), a critically damped oscillator ($\omega_0 = b/2$) and an overdamped oscillator ($\omega_0 < b/2$). The initial conditions are $y(0) = \dot{y}(0) = 0$.

**Proposition 2**. $H_1(\omega) = \dfrac{1}{-\omega^2 + jb\omega + \omega_0^2}$ (12)

$H_2(\omega_1, \omega_2) = -\varepsilon \dfrac{H_1(\omega_1)H_1(\omega_2)}{-(\omega_1 + \omega_2)^2 + jb(\omega_1 + \omega_2) + \omega_0^2}$ (13)

*Proof*:

The rationale behind the computation of these kernels is to work perturbatively as in [10]. The solution has the form:

$y(t) = y_1(t) + y_2(t) + \ldots$ (14)

$y_1(t)$ is the leading term and the solution to the linear part of (11). I.e.

$\ddot{y}_1 + b\dot{y}_1 + \omega_0^2 y_1 = x(t)\theta(t)$ (15)

Thus:

$Y_1(\omega) = H_1(\omega)X(\omega)$ (16)

where $H_1(\omega)$ is given by (12).

$y_2(t)$ is the next correction and is assumed to be small compared to $y_1(t)$. It is the solution to:

$\ddot{y}_2 + b\dot{y}_2 + \omega_0^2 y_2 = -\varepsilon y_1^2$, (17)

$$(-\omega^2 + j\omega b + \omega_0^2)Y_2(\omega) = -\varepsilon \int_{-\infty}^{\infty} \left\{ e^{-j\omega t} \frac{1}{2\pi} \int_{-\infty}^{\infty} e^{j\omega_1 t} Y_1(\omega_1) d\omega_1 \frac{1}{2\pi} \int_{-\infty}^{\infty} e^{j\omega_2 t} Y_1(\omega_2) d\omega_2 \right\} dt$$

$$= -\frac{\varepsilon}{2\pi} \int\!\!\int_{-\infty}^{\infty} \left( \frac{1}{2\pi} \int_{-\infty}^{\infty} e^{j(\omega_1 + \omega_2 - \omega)t} dt \right) Y_1(\omega_1) Y_1(\omega_2) d\omega_1 d\omega_2 =$$

$$= -\frac{\varepsilon}{2\pi} \int\!\!\int \delta(\omega_1 + \omega_2 - \omega) H_1(\omega_1) X(\omega_1) H_1(\omega_2) X(\omega_2) d\omega_1 d\omega_2$$

Thus:

$$Y_2(\omega) = -\frac{\varepsilon}{2\pi} \int\!\!\int \frac{H_1(\omega_1) H_1(\omega_2)}{(-\omega^2 + j\omega b + \omega_0^2)} \delta(\omega_1 + \omega_2 - \omega) X(\omega_1) X(\omega_2) d\omega_1 d\omega_2 \tag{18}$$

Since $Y_2(\omega) = \frac{1}{2\pi} \int\!\!\int H_2(\omega_1, \omega_2) \delta(\omega_1 + \omega_2 - \omega) X(\omega_1) X(\omega_2) d\omega_1 d\omega_2$

then, it follows that

$$H_2(\omega_1, \omega_2) = -\varepsilon \frac{H_1(\omega_1) H_1(\omega_2)}{-(\omega_1 + \omega_2)^2 + jb(\omega_1 + \omega_2) + \omega_0^2}. \tag{19}$$

The degree one kernel $H_1(\omega)$ can further be represented as:

$$H_1(\omega) = \frac{1}{-(\omega - \omega_+)(\omega - \omega_-)} \tag{20}$$

where the poles are

$$\omega_\pm = \frac{jb}{2} \pm \sqrt{\omega_0^2 - \frac{b^2}{4}}.$$

This enables the $a_k^{[1]}$ to be immediately computed from (8) by e.g. the Residue Theorem. The integration contour includes the real axis and the upper infinite semicircle. Thus since both poles $\omega_\pm$ are in the upper half plane, they are inside the contour. Using the Residue Theorem gives:

$$a_k^{[1]} = \frac{\pi j}{\omega_M} \frac{e^{j\tilde{k}\omega_-} - e^{j\tilde{k}\omega_+}}{\omega_+ - \omega_-} \tag{21}$$

For an underdamped oscillator (21) simplifies to:

$$a_k^{[1]} = \frac{\pi}{\omega_M} \frac{e^{-\tilde{k}b/2} \sin(\tilde{k}\omega_R)}{\omega_R} \tag{22}$$

where

$$\omega_R = \sqrt{\omega_0^2 - \frac{b^2}{4}}. \tag{23}$$

For an overdamped oscillator (21) gives:

$$a_k^{[1]} = \frac{\pi}{\omega_M} \frac{e^{-\tilde{k}b/2} \sinh(\tilde{k}\lambda)}{\lambda} \tag{24}$$

where

$$\lambda = \sqrt{\frac{b^2}{4} - \omega_0^2}. \tag{25}$$

For a critically damped oscillator ($\omega_+ = \omega_- = \frac{jb}{2}$)

$$a_k^{[1]} = \frac{\pi}{\omega_M} \tilde{k} e^{-\tilde{k}b/2} \tag{26}$$

The computation of $a_{kl}^{[2]}$ from (9) is rather tedious. It involves applying the Residue Theorem twice and gives:

$$a_{kl}^{[2]} = -\frac{\pi^2 \varepsilon}{\omega_M^2 (\omega_+ - \omega_-)^2} \left[ \frac{e^{j\omega_+(\tilde{k}+\tilde{l})}(1 - e^{-j(2\omega_+ - \omega_-)\tilde{m}})}{\omega_+(2\omega_+ - \omega_-)} + \frac{e^{j\omega_-(\tilde{k}+\tilde{l})}(1 - e^{-j(2\omega_- - \omega_+)\tilde{m}})}{\omega_-(2\omega_- - \omega_+)} \right.$$
$$\left. - \frac{e^{j\omega_-\tilde{k}+j\omega_+\tilde{l}} + e^{j\omega_+\tilde{k}+j\omega_-\tilde{l}} - e^{j\omega_+\tilde{n}} - e^{j\omega_-\tilde{n}}}{\omega_+ \omega_-} \right] \tag{27}$$

where $m = \min(k,l)$, $n = \max(k,l)$.

For an underdamped oscillator (27) gives:

$$a_{kl}^{[2]} = -\frac{\pi^2 \varepsilon e^{-(\tilde{k}+\tilde{l})b/2}}{2\omega_M^2 \omega_0^2 \omega_R^2} \{ \frac{(3\omega_0^2 - b^2)[\cos((\tilde{k}+\tilde{l})\omega_R) - e^{b\tilde{m}/2}\cos((\tilde{k}+\tilde{l}-3\tilde{m})\omega_R)]}{9\omega_0^2 - 2b^2}$$

$$+ \frac{2b\omega_R[\sin((\tilde{k}+\tilde{l})\omega_R) - e^{b\tilde{m}/2}\sin((\tilde{k}+\tilde{l}-3\tilde{m})\omega_R)]}{9\omega_0^2 - 2b^2} \qquad (28)$$

$$+ \cos((\tilde{k}-\tilde{l})\omega_R) - e^{b(\tilde{k}+\tilde{l}-\tilde{n})/2} \cos\omega_R\tilde{n} \}.$$

For a critically damped oscillator:

$$a_{kl}^{[2]} = -\frac{4\pi^2\varepsilon}{\omega_M^2 b^4} \times$$

$$\times \left( e^{-b(\tilde{k}+\tilde{l})/2}(b^2\tilde{k}\tilde{l} + 4b(\tilde{k}+\tilde{l}) + 24) + e^{-b\tilde{m}/2}(b^2\tilde{k}\tilde{l} - 4b(\tilde{k}+\tilde{l}) - b^2\tilde{m}^2 + 12b\tilde{m} - 24) \right)$$

(29)

For an overdamped oscillator, the poles are purely imaginary:

$\omega_\pm = j\lambda_\pm$ where $\lambda_\pm$ are real and positive and are equal to:

$$\lambda_\pm = \frac{b}{2} \pm \lambda$$

where $\lambda$ is given in (25).

From (27):

$$a_{kl}^{[2]} = -\frac{\pi^2\varepsilon}{4\omega_M^2 \lambda^2} [ \frac{e^{-\lambda_+(\tilde{k}+\tilde{l})}(1-e^{(2\lambda_+ - \lambda_-)\tilde{m}})}{\lambda_+(2\lambda_+ - \lambda_-)} + \frac{e^{-\lambda_-(\tilde{k}+\tilde{l})}(1-e^{(2\lambda_- - \lambda_+)\tilde{m}})}{\lambda_-(2\lambda_- - \lambda_+)}$$

$$- \frac{e^{-\lambda_-\tilde{k} - \lambda_+\tilde{l}} + e^{-\lambda_+\tilde{k} - \lambda_-\tilde{l}} - e^{-\lambda_+\tilde{n}} - e^{-\lambda_-\tilde{n}}}{\lambda_+ \lambda_-}] \qquad (30)$$

Note that (30) does not have any singularities even if (possibly) $2\lambda_- = \lambda_+$, since

$$\frac{1 - e^{(2\lambda_- - \lambda_+)\tilde{m}}}{2\lambda_- - \lambda_+} = -\tilde{m} - (2\lambda_- - \lambda_+)\frac{\tilde{m}^2}{2!} - (2\lambda_- - \lambda_+)^2 \frac{\tilde{m}^3}{3!} - \dots.$$

Fig 1. shows the solution to (11) with an input $x(t) = \sin(0.5t)$ and with the parameter values in (11) set as follows: $\varepsilon = 1$, $b = 0.3$ (underdamped oscillator) and $\omega_0 = 2$. The solution is then determined from (3), (22) and (28), where the order of the perturbative terms taken is $n = 2$. Fig. 2 shows the solution to (11) with an input

$x(t) = \theta(t)$ and in this case $\varepsilon = 1$, $b = 5$ (overdamped oscillator) and $\omega_0 = 2$. The solution is determined using (3), (24) and (30). It is evident in both of these cases that the first and second-order approximations are very close to the exact output.

Now consider the underdamped case again where the behaviour is oscillatory. The perturbative approach for the determination of the kernels (thereby permitting the subsequent evaluation of the 'scattering' coefficients) is applicable only for inputs that do not have poles at zero or close to zero in the frequency domain. E.g. in Fig. 1 when the input is $\sin(0.5t)\theta(t)$, the perturbative approach is satisfactory as the input when represented in the frequency domain is $\dfrac{4}{-\omega^2 + 0.25}$ and therefore has poles at $\pm 0.5j$. (In Fig. 2 the system is overdamped and hence non-oscillatory and again the perturbative approach is satisfactory.) However, in the case of a step input to the underdamped system, it will be shown in what follows that use of the perturbative scheme results in the appearance of resonance terms. Thus a modification to the scheme will be proposed and results will confirm its validity.

Consider a step input $K\theta(t)$, where $K$ is a constant. $y_1(t)$ is the leading term in (14) and the solution to:

$$\ddot{y}_1 + b\dot{y}_1 + \omega_0^2 y_1 = K\theta(t) \tag{31}$$

$y_2(t)$ is the next correction and the solution to (17), etc.

The solution of (31) is:

$$y_1(t) = \frac{K}{\omega_0^2}\left(1 - e^{-bt/2}\left(\cos\omega_R t + \frac{b}{2\omega_R}\sin\omega_R t\right)\right)\theta(t) \tag{32}$$

The substitution of (32) into (17) however results in resonance terms in the solution for $y_2(t)$ that are proportional to $t\cos\omega_R t$ and $t\sin\omega_R t$. These terms arise in the

solution due to the fact that $y_1^2$ contains the term

$$-2\frac{K^2}{\omega_0^4}e^{-bt/2}\left(\cos\omega_R t + \frac{b}{2\omega_R}\sin\omega_R t\right)\theta(t)$$

describing oscillations with a resonant frequency $\omega_R$ given in (23). The problem with the $t\cos\omega_R t$ and $t\sin\omega_R t$ terms is that their presence violates, at large values of $t$, the assumption that the second term $y_2(t)$ is small in comparison to the leading term, $y_1(t)$. It should also be noted that the failure of the perturbative scheme for the underdamped case with a step input is a direct consequence of the implicit and incorrect assumption in the scheme that the frequency of oscillation of the anharmonic oscillator is the same as that of the damped harmonic oscillator. Indeed, suppose that the frequency of oscillations of the anharmonic oscillator is $\omega_R + \delta\omega$, where $\delta\omega$ is a small correction. The resonance terms in the solution arise from an expansion of the type

$$\cos(\omega_R + \delta\omega)t \approx \cos\omega_R t - t\delta\omega\sin\omega_R t$$

which is obviously not legitimate for large values of $t$ [15]. Although this fact is well known, it is often overlooked [16].

In the frequency domain, the failure of the perturbative scheme can be explained as follows: if the input has a pole at zero, then from (18):

$$Y_2(\omega) = \frac{1}{2\pi}\int_{-\infty}^{\infty} H_2(\omega_1, \omega - \omega_1)X(\omega_1)X(\omega - \omega_1)d\omega_1 \qquad (33)$$

The Residue Theorem is applicable since $H_n(\omega_1,\ldots,\omega_n)$ is zero (or tends fast enough to zero) for large values of the arguments. The contour consists of the real axis - avoiding the real poles by small detours and the infinite upper semicircle. Now there is, in general, a contribution to the integral from the detours around the real poles at $\omega_1 = 0$ and $\omega_1 = \omega$, proportional to $H_2(0,\omega)$ and $H_2(\omega,0)$. But these, however, have double poles at $\omega = \omega_{\pm}$ from (19) resulting in resonant terms and violating the assumption regarding the smallness of $y_2(t)$.

The perturbative scheme for the determination of the Volterra kernels therefore needs to be modified to account for the different frequency of the anharmonic oscillations to those of the damped harmonic oscillator. To this end, equations (31) and (17) are replaced by:

$$\ddot{y}_1 + b\dot{y}_1 + \Omega_0^2 y_1 = K\theta(t) \tag{34}$$

$$\ddot{y}_2 + b\dot{y}_2 + \omega_0^2 y_2 = -\varepsilon y_1^2 + (\Omega_0^2 - \omega_0^2)y_1 \tag{35}$$

where $\Omega_0$ is a constant whose value is determined so as to enforce the condition that no sources with frequency $\omega_R$ appear on the right hand side of (35). The solution to (34) is:

$$y_1(t) = \frac{K}{\Omega_0^2}\left(1 - e^{-bt/2}\left(\cos\Omega_R t + \frac{b}{2\Omega_R}\sin\Omega_R t\right)\right)\theta(t) \tag{36}$$

where

$$\Omega_R = \sqrt{\Omega_0^2 - \frac{b^2}{4}}. \tag{37}$$

The condition for the elimination of sources with a frequency $\omega_R$ on the right hand side of (35) is:

$$\Omega_0^2 - \omega_0^2 = \frac{2\varepsilon K}{\Omega_0^2} \tag{38}$$

This eliminates the terms with the resonance frequency $\omega_R$ from the right-hand side of (35). (38) results in a quadratic equation for $\Omega_0^2$ with a solution:

$$\Omega_0^2 = \frac{1}{2}\left(\omega_0^2 + \sqrt{\omega_0^4 + 8\varepsilon K}\right) \tag{39}$$

which is the first correction to the frequency due to the anharmonicity.

Equations (38) and (39) show that due to the nonlinearity in (11), the frequency of the anharmonic oscillations depends on the input through the amplitude $K$.

The results for a step input are presented on Fig. 3: $x(t) = \theta(t)$, $K = 1$; $\varepsilon = 1$, $b = 0.3$, $\omega_0 = 2$. This case clearly shows that (13) and (28) are not valid. However, it is evident that the first (linear) approximation (22) at the perturbed frequency, $\Omega_0$ (39) matches the frequency of the anharmonic oscillator.

## 5. Numerical algorithm for computation of the 'scattering' coefficients

In order to compute the 'scattering' coefficients numerically, it is necessary to obtain a recursive formula. In this section for simplicity, the second order differential equation in (11) will be considered to illustrate the approach, but obviously the method can be applied to any ordinary differential equation.

As a first step, the expression for the first approximation, $y_1(t)$, and its derivatives have to be substituted into (15).

$$y_1(t) = \sum_{k=0}^{N} a_k^{[1]} x(t - \tilde{k}) \theta(t - \tilde{k}) \tag{40}$$

There are several possibilities for approximating the first derivative. Using

$\dot{f}(t) = \dfrac{f(t+T) - f(t)}{T}$ with a sampling time $T = \tilde{1} = \pi / \omega_M$ results in:

$$\dot{y}_1(t) = \frac{1}{T} \sum_{k=0}^{N} (a_{k+1}^{[1]} - a_k^{[1]}) x(t - \tilde{k}) \theta(t - \tilde{k}) \tag{41}$$

$$\ddot{y}_1(t) = \frac{1}{T^2} \sum_{k=0}^{N} (a_{k+2}^{[1]} - 2a_{k+1}^{[1]} + a_k^{[1]}) x(t - \tilde{k}) \theta(t - \tilde{k}) \tag{42}$$

(All coefficients with an index greater than $N$ or less than zero is assumed to be zero.)

The substitution of (40), (41) and (42) into (15) gives the recurrence equation:

$$a^{[1]}_{k+2} + (bT - 2)a^{[1]}_{k+1} + (1 + \omega_0^2 T^2 - bT)a^{[1]}_k = 0 \tag{43}$$

To confirm the consistency of equation (43) with the exact solution, (21), consider the following:

The characteristic equation for (43) is:

$$z^2 + (bT - 2)z + 1 + \omega_0^2 T^2 - bT = 0$$

with roots

$$z_\pm = 1 + T\left(-\frac{b}{2} \pm \sqrt{\frac{b^2}{4} - \omega_0^2}\right) = 1 + jT\omega_\pm \tag{44}$$

The general solution of (43) is therefore:

$$a^{[1]}_k = c_+ z_+^k + c_- z_-^k$$

for two arbitrary constants $c_\pm$ that may be determined from the initial data. The initial condition $y_1(0)=0$ gives $a^{[1]}_0 = 0$ and therefore $c_+ = -c_- = c$.

In order to determine $a^{[1]}_1$ it is more convenient to use the following approximation for the second derivative:

$$\ddot{y}_1(t) = \frac{\dot{y}_1(t) - \dot{y}_1(t - T)}{T} \tag{45}$$

Thus:

$$\ddot{y}_1(t) = \frac{1}{T^2}\sum_{k=0}^{N}(a^{[1]}_{k+1} - 2a^{[1]}_k + a^{[1]}_{k-1})x(t - \tilde{k})\theta(t - \tilde{k}) \tag{46}$$

Now, equation (15) with zero initial conditions gives:

$$\ddot{y}_1(0) = x(0) \tag{47}$$

Thus from (46) and (47) assuming $a^{[1]}_{-1} = 0$ yields $a^{[1]}_1 = T^2$. Hence,

$$c = \frac{T}{j(\omega_+ - \omega_-)} = \frac{T}{2\sqrt{\frac{b^2}{4} - \omega_0^2}} \tag{48}$$

Furthermore, when $T$ is small, $z_{\pm} = 1 + jT\omega_{\pm} \approx e^{j\tilde{\omega}_{\pm}}$. So finally:

$$a_k^{[1]} = c(z_+^k - z_-^k) \approx \frac{\pi(e^{jk\tilde{\omega}_+} - e^{jk\tilde{\omega}_-})}{j\omega_M(\omega_+ - \omega_-)}$$

which corresponds exactly with (21). (Note $T = \frac{\pi}{\omega_M}$)

Therefore, the 'scattering' coefficients for the first-order approximation may be computed numerically from the recurrence equation

$$\begin{aligned}
a_0^{[1]} &= 0 \\
a_1^{[1]} &= T^2 \\
a_{k+2}^{[1]} &= (2 - bT)a_{k+1}^{[1]} - (1 + \omega_0^2 T^2 - bT)a_k^{[1]}
\end{aligned} \qquad (49)$$

The 'scattering' coefficients $a_{kl}^{[2]}$ for the second-order approximation, $y_2(t)$, may be computed in a similar fashion:

$$y_2(t) = \sum_{l=0}^{N}\sum_{k=0}^{N} a_{kl}^{[2]} x(t-\tilde{k})x(t-\tilde{l})\theta(t-\tilde{k})\theta(t-\tilde{l}) \qquad (50)$$

Substituting:

$$\dot{f}(t) = \frac{f(t) - f(t-T)}{T}$$

for the derivatives gives:

$$\dot{y}_2(t) = \frac{1}{T}\sum_{l=0}^{N}\sum_{k=0}^{N}(a_{kl}^{[2]} - a_{k-1,l-1}^{[2]})x(t-\tilde{k})x(t-\tilde{l})\theta(t-\tilde{k})\theta(t-\tilde{l}) \qquad (51)$$

$$\ddot{y}_2(t) = \frac{1}{T^2}\sum_{l=0}^{N}\sum_{k=0}^{N}(a_{kl}^{[2]} - 2a_{k-1,l-1}^{[2]} + a_{k-2,l-2}^{[2]})x(t-\tilde{k})x(t-\tilde{l})\theta(t-\tilde{k})\theta(t-\tilde{l}) \qquad (52)$$

Substitution of (50)-(52) and the solution (40) for $y_1(t)$ (with (49)) into (17) results in the following recurrence relation for the 'scattering' coefficients $a_{kl}^{[2]}$:

$$a_{kl}^{[2]} = \frac{(2+bT)a_{k-1,l-1}^{[2]} - a_{k-2,l-2}^{[2]} - \varepsilon T^2 a_k^{[1]} a_l^{[1]}}{1 + bT + \omega_0^2 T^2} \qquad (53)$$

where all data with negative indices is assumed zero. The numerical examples (Fig. 4 and 5) confirm that very good accuracy is achieved when the 'scattering' coefficients are computed numerically as compared to the results obtained when the exact formulae (28), (30).

Clearly, a numerical scheme may be developed for any nonlinear differential equation following the procedured outlined above.

### 6. Computation of the 'scattering' coefficients- inverse problem

Suppose that the nonlinear differential equation, describing the system is NOT explicitly given. Instead, suppose that the output $Y(\omega_i)$ and the input $X(\omega_i)$ are given at a discrete number of points $\omega_i$, covering the interval $[0, \omega_M]$, $i = 0,..., N$. This is often the case when the frequency response of an unknown system is investigated, e.g. with a network analyser. The inverse problem is to compute the 'scattering' coefficients from the given input-output data in the frequency domain. For a nonlinear system, this is a nontrivial problem. One attempt in this direction has been reported in [17]. In the present work, an insight will be given into another approach. However, it will be restricted to a linear system.

Suppose that the real or imaginary parts (or both) of the transfer function of a linear system $H_1(\omega)$ are measured at a number of points, $\omega_i$ - i.e. suppose that the following data is known from measurements:

$$F_i^{(1)} = \operatorname{Re} H(\omega_i), \quad i = 1, 2, \ldots, N_1 \tag{54}$$

$$F_i^{(2)} = \operatorname{Im} H(\omega_i), \quad i = 1, 2, \ldots, N_2 \tag{55}$$

Let $a = [a_0 \quad a_1 \quad \ldots \quad a_N]^T$, $M_{ik}^{(1)} = \cos k\tilde{\omega}_i$, $M_{ik}^{(2)} = -\sin k\tilde{\omega}_i$, $k = 1,..,N$. The upper index '[1]' of the 'scattering' coefficients will be omitted for simplicity throughout this section.

Then from (6) and (7):

$$F^{(1)} = M^{(1)}a + E^{(1)} \tag{56}$$

$$F^{(2)} = M^{(2)}a + E^{(2)} \tag{57}$$

$E^{(1,2)}$ represent the errors that arise due to limiting the summation to a finite number of terms, $N$. Note that (56) represents a formulation for the least-squares method and therefore provides the best approximation for $a$ by minimizing the error $E^{(1)T}E^{(1)}$:

$$a = \left(M^{(1)T}M^{(1)}\right)^{-1}M^{(1)T}F^{(1)} \tag{58}$$

The system in (57), may be used to retrieve all the elements of $a$ with the exception of $a_0$ since $M_{i0}^{(2)} \equiv 0$. However, if data for both the real and the imaginary parts is available, then (56) and (57) may be merged to yield:

$$F = Ma + E \tag{59}$$

with

$$F = \begin{bmatrix} F^{(1)} \\ F^{(2)} \end{bmatrix}, \quad M = \begin{bmatrix} M^{(1)} \\ M^{(2)} \end{bmatrix}, \quad E = \begin{bmatrix} E^{(1)} \\ E^{(2)} \end{bmatrix}. \tag{60}$$

and the minimal error, $E^T E$, for (59) is achieved with:

$$a = \left(M^T M\right)^{-1} M^T F \tag{61}$$

In order to apply the least squares method with (59) it is essential that:

$$N < \dim(F) = N_1 + N_2,$$

$$\tilde{N} \geq T_M$$

where $[0,T_M]$ is the time-interval under consideration. Thus, if more measurement data that is available, then more coefficients $\{a_k\}$ may be found and consequently, a better resultant interpolation of (5) may be achieved.

Once the vector of coefficients $a$ is found from either (58) or (61), the response of the system for an arbitrary input from (3) is:

$$y(t) = \sum_{k=0}^{N} a_k x(t - \tilde{k}) \theta(t - \tilde{k}) \quad (62)$$

The precision of (62) clearly increases with increasing $N$. Thus, the accuracy of the method is only limited by the volume of available measurement data. Note that the vector of 'scattering' coefficients, $a$, is obtained from the data vector, $F$, through a simple linear transformation (61)! Furthermore, because of the use of the Fourier Series as opposed to a pole-residue representation, the problems with ill-conditioning do not arise. The application of this approach to RF interconnect networks will be studied further elsewhere [18].

## 7. Conclusions

The input-output behaviour of a nonlinear system in both the frequency and time domain may be determined with the aid of the 'scattering' coefficients. The response in both domains can be represented in a closed form for an arbitrary input. Algorithms are outlined for both the exact and numerical computation of the 'scattering' coefficients for the specific case of the anharmonic oscillator. It is shown that the perturbative approach for the determination of the Volterra kernels is not valid for the case when the anharmonic oscillator system exhibits oscillatory behaviour and the input has a pole at the origin (in the frequency domain). For this case, the frequency of the oscillations depends on the amplitude of the input. Consequently, the entire perturbative scheme, that is employed to determine the Volterra Series kernels needs to be appropriately modified in order to incorporate this dependence. Finally, an approach for the determination of the 'scattering' coefficients from measured data is given. The procedure is explained relative to the linear case and is presented as a starting point for future work.

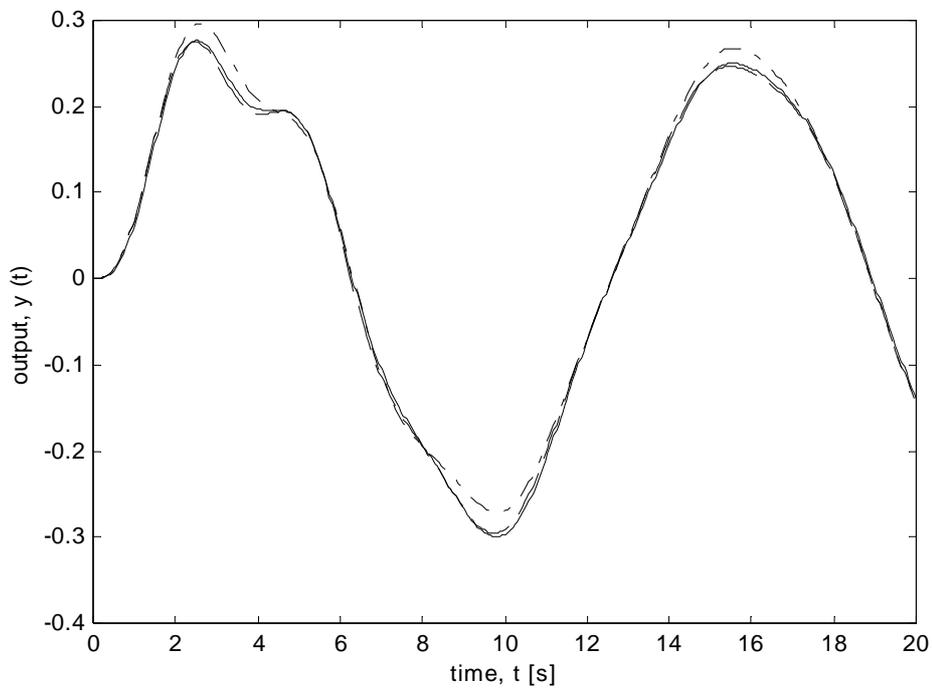

Fig. 1. Underdamped anharmonic oscillator with input sin(0.5$t$). Exact solution (solid line), First-order (linear) approximation $y_1(t)$ (dash-dotted line); Second-order approximation $y_1(t)+y_2(t)$ (dashed line).

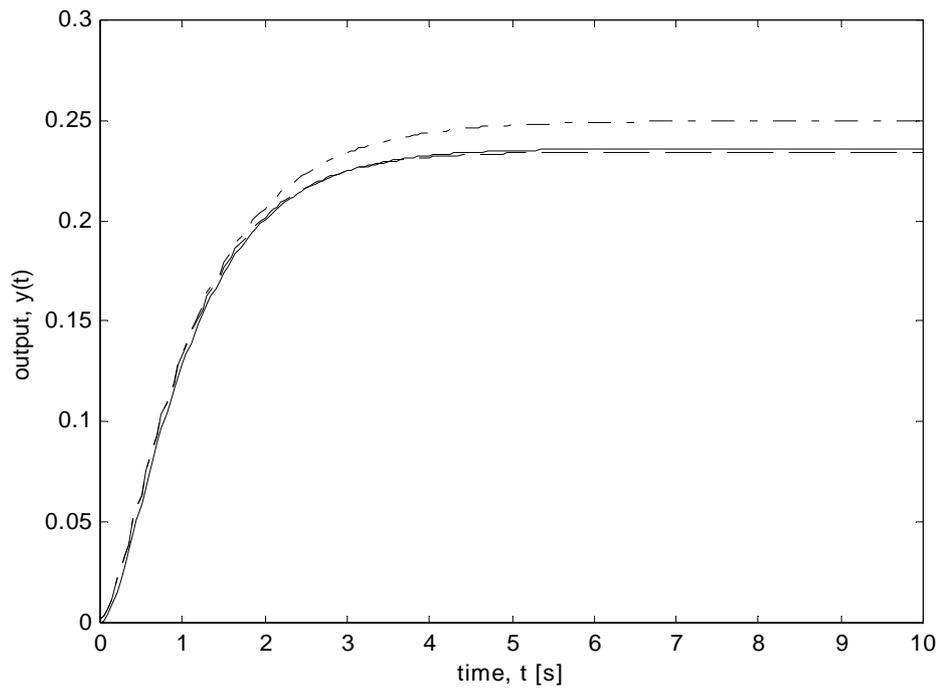

Fig. 2. Overdamped anharmonic oscillator-step input: Exact solution (solid line), First-order (linear) approximation $y_1(t)$ (dash-dotted line); Second-order approximation $y_1(t) + y_2(t)$ (dashed line).

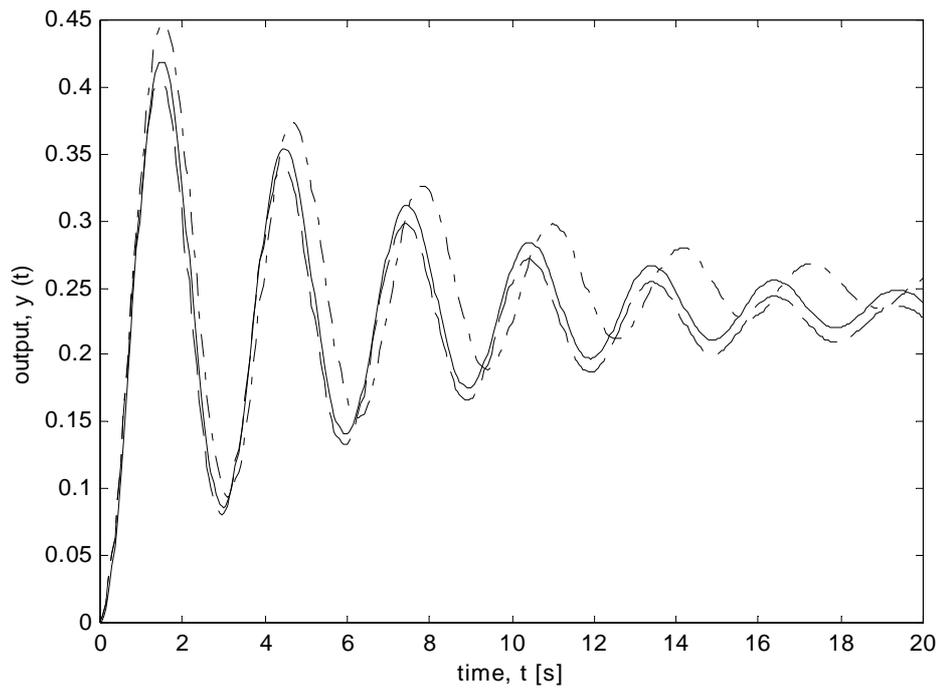

Fig. 3. Underdamped anharmonic oscillator- step input: Exact solution (solid line), First-order(linear) approximation $y_1(t)$ with nonperturbed frequency (23) (dash-dotted line); First-order (linear) approximation $y_1(t)$ with perturbed frequency (37) (dashed line).

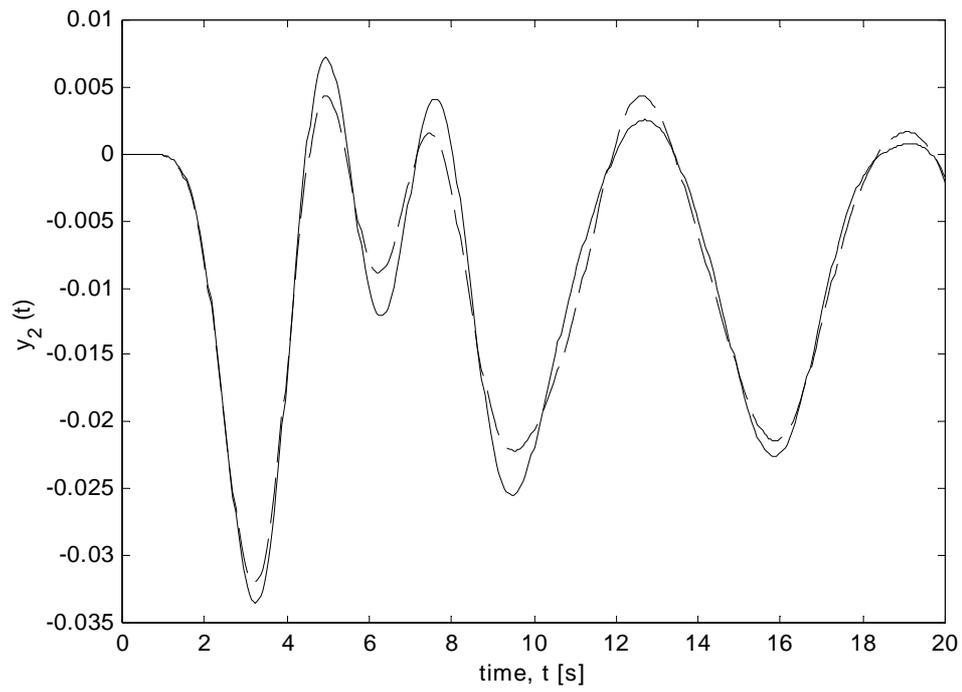

Fig. 4. The second-order approximation $y_2(t)$: underdamped oscillator, sine input.
Exact expression (28) (solid line), numerical approximation (53) (dashed line).

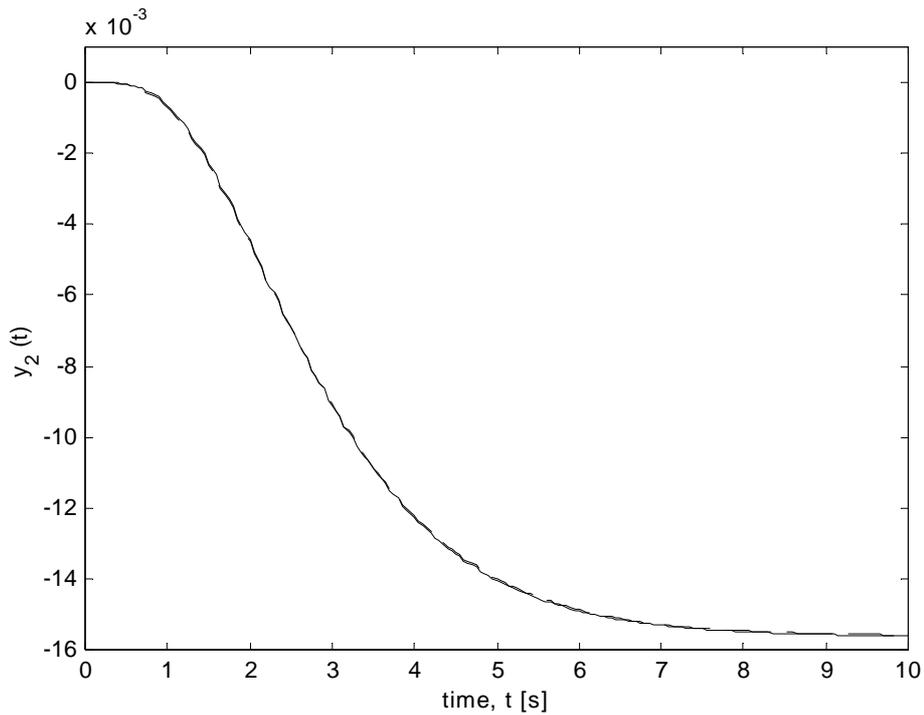

Fig. 5. The second-order approximation $y_2(t)$: overdamped oscillator, step input. Exact expression (30) (solid line), numerical approximation (53) (dashed line).

**REFERENCES**

[1] B. Ulriksson, Conversion of frequency-domain data to time-domain, *Proceedings of the IEEE* 74, no.1 (1986) 74-77.

[2] T. Brazil, Causal-convolution—a new method for the transient analysis of linear Systems at microwave frequencies, *IEEE Transactions on microwave theory and techniques* 43, no.2 (1995) 315-323.

[3] M. Condon and E. Dautbegovic, "Efficient modelling of interconnects in high-speed circuits", *Proc. European Conference on Circuit theory and Design*, Krakow, (2003).

[4] J.D. Jackson, *Classical Electrodynamics*, 2$^{nd}$ ed., Wiley, New York, 1975.

[5] J.G. McDaniel, Applications of the causality condition to acoustic reflections, *Proceedins of 1997 ASME Design Engineering Technical Conferences*, Sacramento, California, DETC97/VIB-4133, Seprember 14-17, 1997.

[6] V. Volterra, *Theory of functionals and of integrodifferential equations,* Dower, New York, 1959.

[7] M. Schetzen, *The Volterra and Wiener Theories of Nonlinear Systems,* Wiley, New York, 1980.

[8] A.I. Mees, The describing Function Matrix, *Journal of the Institute of Mathematics and its Appications*, 10 (1972) 49-67.

[9] J.C. Peyton Jones and S.A. Billings, Describing functions, Volterra series and the analysis of nonlinear systems in the frequency domain, *International Journal of Control* 53, no. 4 (1991) 871-887.

[10] S.A. Billings and J.C. Peyton Jones, Mapping nonlinear integro-differential equations into the frequency domain, *International Journal of Control* 52, no. 4 (1990) 863-879.

[11] J.C. Peyton Jones and S.A. Billings, Interpretation of non-linear frequency response functions, *International Journal of Control*, 52, no. 2 (1990) 319-346.

[12] J.C. Peyton Jones and S.A. Billings, Recursive algorithm for computing the frequency response of a class of non-linear difference equations models, *International Journal of Control* 50, no. 5 (1989) 1925-1940.

[13] F. Bassani and S. Scandolo, Dispersions and sum rules in nonlinear optics, *Phys. Rev. B* 44, no.16 (1991) 8446-8453.

[14] F. Smet and A. van Groenendael, Dispersion relations for *n*-order nonlinear phenomena, *Phys. Rev. A* 19, no.1 (1979) 334-337.

[15] L.D. Landau and E.M. Lifshitz, *Mechanics*, Vol.1 of *Course of Theoretical Physics*, Pergamon Press, Oxford, 1976.

[16] S. Mandal "Classical damped quartic anharmonic oscillator: a simple analytical approach", *Int. J. of Non-linear Mechanics* 38 (2003) 1095-1101.

[17]. Verspecht, D. Schreurs, A. Barel and B. Nauwelaers, Black Box Modelling Of Hard Nonlinear Behavior In The Frequency Domain, *IEEE MTT-S Digest*, San Francisco, USA, pp. 1735-1738, 17-21 June 1996.

[18]. M. Condon, R. Ivanov and C. Brennan, A causal model for linear RF systems developed from frequency domain measured data, IEEE TCAS II (2005), to appear.